\documentclass[12pt]{article}
\usepackage{latexsym, amssymb, amsmath, amscd, amsfonts,mathrsfs}

\newtheorem{theorem}{Theorem}[section]
\newtheorem{lemma}[theorem]{Lemma}

\numberwithin{equation}{section}

\def\qed{\hfill \rule{4pt}{7pt}}

\def\pf{\noindent {\it Proof.} }

\begin{document}
\begin{center}
{\large {\bf On Han's Hook Length Formulas for Trees}}

\vskip 6mm

{\small William Y.C. Chen$^1$, Oliver X.Q. Gao$^2$
and Peter L. Guo$^3$\\[%
2mm] Center for Combinatorics, LPMC-TJKLC\\
Nankai University, Tianjin 300071,
P.R. China \\[3mm]
$^1$chen@nankai.edu.cn, $^2$oliver@cfc.nankai.edu.cn, $^3$lguo@cfc.nankai.edu.cn \\[0pt%
] }
\end{center}

\begin{abstract}
Recently, Han obtained two hook length
 formulas for binary trees
and asked for   combinatorial proofs.
 One of Han's formulas
has been generalized to $k$-ary trees by Yang.
Sagan has found a
probabilistic proof of Yang's extension.
We give combinatorial
proofs of  Yang's formula for $k$-ary trees and the other
formula of Han for binary trees. Our bijections are based on
the structure of  $k$-ary trees with staircase labelings.
\end{abstract}
\vskip 3mm



\section{Introduction}

Motivated by the hook length formula
 of Postnikov \cite{Pos}, Han
\cite{Han2}  discovered two hook length formulas for binary trees.
Han's proofs are based on
recurrence relations. He raised the
question of finding combinatorial proofs of these two
formulas \cite{Han1,
Han2}. Yang \cite{Yang}
 generalized one of Han's formulas to $k$-ary
trees by using  generating functions. A probabilistic
proof of Yang's formula  has been found
by Sagan \cite{Sag}. By
extending Han's expansion technique
 to $k$-ary trees, Chen, Gao
and Guo \cite{Che1} gave another
proof for Yang's formula. The
objective of this paper
is to give combinatorial proofs of Yang's
 formula  for $k$-ary trees and  the other
  formula of Han
for binary trees.

Recall that a
$k$-ary tree is a rooted unlabeled tree where each vertex has
exactly $k$ subtrees in linear order, where we allow a subtree to be
empty. When $k=2$ (resp., $k=3$), a $k$-ary tree is called  a binary
(resp., ternary) tree. A complete  $k$-ary tree is a $k$-ary tree
for which each internal vertex has
exactly $k$ nonempty subtrees.
The hook length of a vertex $u$ in a $k$-ary tree $T$,
denoted by
$h_u$, is the number of vertices
 of the subtree rooted at $u$. The hook length
multi-set $\mathcal {H}(T)$ of $T$ is
 defined to be the multi-set of
hook lengths of all vertices of $T$. For
 example, Figure \ref{h}
gives an illustration of the hook length multi-set of a binary tree.
\begin{figure}[h,t]
\setlength{\unitlength}{0.5mm}
\begin{center}
\begin{picture}(110,50)
\put(-10,15){\circle*{2}}\put(0,30){\circle*{2}}\put(10,45){\circle*{2}}
\put(10,15){\circle*{2}}\put(20,30){\circle*{2}}\put(30,15){\circle*{2}}
\put(-10,15){\line(2,3){10}}\put(0,30){\line(2,3){10}}
\put(10,45){\line(2,-3){10}}
\put(10,15){\line(2,3){10}}\put(20,30){\line(2,-3){10}}
\put(6,0){\shortstack{$T$}}\put(50,25){$\mathcal
{H}(T)=\{1,1,1,2,3,6\}$}
\end{picture}\caption{The multi-set of
hook lengths of a binary tree.}
\end{center}  \label{h}
\end{figure}

Let $B_n$ be the set of all binary trees with $n$ vertices. Han
\cite{Han2} discovered the following  formulas.
He also gave derivations of
these formulas in \cite{Han1} by using
 the expansion technique.

\begin{theorem}[\mdseries{Han \cite{Han2}}]\label{H}
For each positive integer $n$, we have
\begin{align}\label{1}
\sum_{T\in {B}_n} \frac{1}{\prod_{h\in
 \mathcal{H}(T)}h
2^{h-1}}=\frac{1}{n!}
\end{align} and
\begin{align}\label{2}
\sum_{T\in {B}_n} \frac{1}{\prod_{h\in
 \mathcal{H}(T)}(2h+1)
2^{2h-1}}=\frac{1}{(2n+1)!}.
\end{align}
\end{theorem}
As pointed out by Han \cite{Han2},  the above
two formulas have a special feature
that the hook lengths appear as  exponents. Yang
\cite{Yang} extended the
above formula (\ref{1}) to $k$-ary trees.

\begin{theorem}[\mdseries{Yang \cite{Yang}}]\label{Y}
For any positive integers $n$ and $k$, we have
\begin{align}\label{4}
\sum_{T}\prod\limits_{h\in
 \mathcal{H}(T)}\frac{1}{hk^{h-1}}=\frac{1}{n!},
\end{align}
where the sum ranges over  $k$-ary trees with $n$ vertices.
\end{theorem}

To give a combinatorial proof of (\ref{4}),
we shall define a set
$S(n,k)$ of staircase arrays on $[k]=\{1,2,\ldots, k\}$. More
precisely, we shall represent an array in
$S(n,k)$ in the form $(C_0,C_1,\ldots, C_{n-1})$,
where $C_0=\emptyset$ and for $1\leq i \leq n-1$,
 $C_i$
is a vector of length $i$ with
each entry in $[k]$.

We shall show
 that the sequences in $S(n,k)$ are in
 one-to-one correspondence to
 $k$-ary trees with $n$  vertices
 whose labels form a staircase sequence in
 $S(n,k)$. Such  $k$-ary trees  are called
   $k$-ary trees with   staircase labelings.
   This  leads
to a bijective proof of formula (\ref{4}).
Based on this bijection,
we further obtain
 a combinatorial interpretation of formula
(\ref{2}).

\section{A combinatorial proof of (\ref{4})}

Our combinatorial proof of Yang's formula \eqref{4} is based on the
following reformulation
\begin{equation}\label{10}
\sum_{T}\frac{n!k^{1+2+\cdots+n}}{\prod_{h\in
 \mathcal{H}(T)}hk^h}=k^{1+2+\cdots+(n-1)}.
\end{equation}
It is clear that the right-hand side of (\ref{10}) equals the number
of sequences in $S(n,k)$. As will be seen,
 the left hand-side of
\eqref{10} equals the number of
 $k$-ary trees with $n$ vertices
 whose labels form a
staircase array. Such a $k$-ary tree is called a $k$-ary tree with a
staircase labeling. Furthermore, we shall give a one-to-one
correspondence between $S(n,k)$ and the set of $k$-ary
trees with $n$ vertices associated with
staircase labelings.

More precisely, a staircase labeling of a $k$-ary tree is defined as
follows: The labels are vectors on $[k]$
 with distinct lengths. In
other words, the labels can be written as
$C_0,C_1,\ldots,C_{n-1}$,
where $C_i$ is a vector on $[k]$
of length $i$. Moreover, we
impose the following restrictions:
for any vertex $u$ with label
$C_i$  and a descent (not necessarily a child)
 $v$
 with label $C_j$, we have $i< j$,
 that is, the labels on any path from the root to a leaf
 have increasing lengths;
 and the $(i+1)$-st entry of  $C_j$ is determined by
 the relative position of the  child of
 $u$  on the path from $u$ to $v$ among its siblings.
 To be more specific, if the $r$-th child of $u$ is on the
 path from $u$ to $v$, then the $(i+1)$-st entry is
 set to be $r$.

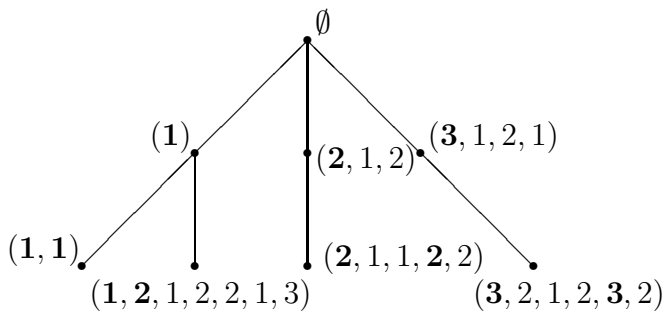
\begin{figure}[h,t]
\setlength{\unitlength}{1.0mm}
\begin{center}
\begin{picture}(90,42)
 \put(15,5){\line(1,1){15}}\put(15,5){\circle*{1}}
 \put(30,5){\line(0,1){15}}\put(30,5){\circle*{1}}
 \put(45,5){\line(0,1){15}}\put(45,5){\circle*{1}}
 \put(75,5){\line(-1,1){15}}\put(75,5){\circle*{1}}
 \put(30,20){\line(1,1){15}}\put(30,20){\circle*{1}}
 \put(45,20){\line(0,1){15}}\put(45,20){\circle*{1}}
 \put(60,20){\line(-1,1){15}}\put(60,20){\circle*{1}}
 \put(45,35){\circle*{1}}
\put(46,36){$\emptyset$} \put(5,6){$(\textbf{1},\textbf{1})$}
\put(16,0){$(\textbf{1},\textbf{2},1,2,2,1,3)$}
\put(47,5){$(\textbf{2},1,1,\textbf{2},2)$}
\put(67,0){$(\textbf{3},2,1,2,\textbf{3},2)$}\put(24,21){$(\textbf{1})$}
\put(46,18){$(\textbf{2},1,2)$} \put(61,21){$(\textbf{3},1,2,1)$}
\end{picture}\caption{A staircase labeling of a
ternary tree}\label{t}
\end{center}
\end{figure}

 For example, Figure \ref{t}
 gives a staircase labeling of a
ternary tree, where the label of any vertex,
 the entries  that
are determined by the labels of its ancestors are
written in boldface.

Let $I(n,k)$ denote the set of  $k$-ary trees
with $n$ vertices associated with staircase labelings.
The following lemma shows that $|I(n,k)|$ is
equal to the left-hand side of (\ref{10}).

\begin{lemma}\label{bb}
For $n\geq 1$,
\begin{equation}\label{b}
|{I}(n,k)|=\sum_{T}\frac{n!k^{1+2+\cdots+n}}{\prod_{h\in
 \mathcal{H}(T)}hk^h},
\end{equation}
where the sum ranges over $k$-ary trees with $n$ vertices.
\end{lemma}

\pf Let $P\in I(n,k)$ be a $k$-ary tree with a staircase labeling.
Suppose that the labels of $P$ are $C_0, C_1, \ldots, C_{n-1}$,
where $C_i$ is a vector of length $i$. Define $Q$ to be the $k$-ary
tree obtained from $P$ by replacing a label $C_i$ with $i$. Clearly,
$Q$ is an increasing $k$-ary tree in the sense that the label of any
internal vertex is smaller than the labels of its children.

We shall consider the question of determining the number of $k$-ary
trees $P$ in $I(n,k)$ that correspond to a given
increasing $k$-ary
tree $Q$. Clearly, $P$ and $Q$ have the same underlying $k$-ary
tree, denoted by $T$. In other words, we shall compute the
number of staircase labelings of a $k$-ary tree $T$
with given label length for each vertex.
For any vertex $u$ of $T$, let $f_u$ denote the
number of vertices on the path from the root to $u$.
We claim that
there are
\begin{equation}\label{a}
k^{1+\cdots+n-\sum_{u\in T}f_u}
\end{equation}
 staircase labelings of $T$ such that the label corresponding
 to a vertex with label $i$ in $Q$ has length $i$.
 To prove (\ref{a}),
 let $u_i$ be the vertex of
$Q$ with label $i$.
Assume that $P$ is a $k$-ary tree with a
staircase labeling for which the vertex $u_i$
has label $C_i$. Notice that $C_i$ is of length $i$.
Recalling the definition of a staircase labeling,
we need to determine how many entries in $C_i$ that are
determined by the ancestors of $u_i$.
It can be seen that there are  $f_u-1$ entries of
$C_i$ that are determined by the ancestors
of $u_i$. The other entries can be any choice of the
elements in $[k]$. Consequently,
there are $k^{i+1-f_u}$
choices for $C_i$. This implies \eqref{a}.

Note that the number  in \eqref{a} does not
depend on the specific increasing labeling of
the $k$-ary tree $T$. To compute the
number of $k$-ary trees
with staircase labelings, it suffices to
determine the number of increasing labelings of
$T$. It is known
that the number of increasing
labelings of $T$ equals
 \[ \frac{n!}{\prod_{h\in \mathcal{H}(T)}h}\]
   see, for example, Gessel and Seo \cite{Gel}.
So we deduce that
\begin{equation}\label{8}
|{I}(n,k)|=\sum_{T}\frac{n!}{\prod_{h\in
 \mathcal{H}(T)} h}k^{1+\cdots+n-\sum_{u\in T}f_u},
\end{equation}
where $T$ ranges over $k$-ary trees with $n$ vertices.

We now need to establish the following relation
\begin{equation}\label{9}
\sum_{u\in T}h_u=\sum_{u\in T} f_u.
\end{equation} This can be justified
 by observing that both sides of
\eqref{9} count the number of  ordered  pairs $(u,v)$, where
$v$ is a descendant of $u$ in $T$ and  we adopt the assumption  that $u$ is a
descendant of itself. Substituting
(\ref{9}) into (\ref{8}), we
arrive at \eqref{b}. This completes the proof. \qed

We have the following correspondence.

\begin{theorem}\label{th}
There is a bijection between $ S(n,k)$ and ${I}(n,k)$.
\end{theorem}

\pf The map $\varphi$ from ${I}(n,k)$ to $S(n,k)$ is
straightforward, that is,  for    $P \in{I}(n,k)$
 with a labeling set
$ \{C_0,C_1,\ldots, C_{n-1} \}$, define
\[\varphi(P)=(C_0,C_1,\ldots, C_{n-1}).\]

We now proceed  to give the  reverse map
$\phi$ from $ S(n,k)$ to
${I}(n,k)$. Given a sequence
 $(C_0,C_1,\ldots, C_{n-1})$ in
$S(n,k)$, we aim to construct a $k$-ary tree
with $n$ vertices
associated with a staircase labeling   $\{C_0,C_1,\ldots,
C_{n-1}\}$.

The map $\phi$ can be described as a recursive procedure.
 Let $v_0$ be a
 vertex with label
$C_0=\emptyset$. Clearly, $v_0$
and its label $C_0$ form a $k$-ary
tree with a staircase labeling.
Let $C_1=(c_1)$. Adding a vertex $v_1$ as
the $c_1$-th child of $v_0$ and assigning
the label $C_1$ to $v_1$, we get
 a $k$-ary tree labeled by $C_0$ and $C_1$, denoted by
$P_1$. It can be easily checked that $P_1$ is a $k$-ary
tree with a
staircase labeling.
Assume that $P_{m-1}$ ($m\geq 2$) is a $k$-ary
tree with a staircase labeling with
vertices $v_0,v_1,\ldots,v_{m-1}$ such that
 for $0\leq i\leq m-1$ the vertex
  $v_i$ has label $C_i$.  Now we construct
 a $k$-ary tree with a staircase labeling,
 denoted by $P_{m}$,
 by adding the vertex $v_{m}$ to $P_{m-1}$ and assigning the
 label $C_m$ to $v_m$.

To determine the position of $v_m$, we start at
the root $v_0$. Let $C_m=(c_1,c_2,\ldots,c_m)$.
If  the $c_{1}$-th child of $v_0$ is empty,
then we add
the  vertex $v_m$
to $P_{m-1}$ as the $c_{1}$-th child of $v_0$.
Otherwise,
we arrive at the $c_{1}$-th child of
$v_0$, denoted by $v_{i_0}$. Consider the
the label $C_{i_0}$ of $v_{i_0}$.
If  the $c_{i_0+1}$-th child of $v_{i_0}$ is empty,
then we add
the  vertex $v_m$
to $P_{m-1}$ as the $c_{i_0+1}$-th child of $v_{i_0}$.
Otherwise,
we arrive at the $c_{i_0+1}$-th child of
$v_{i_0}$. Repeating this process, we finally
get a  $k$-ary tree $P_m$
with labeled by $C_0,C_1,\ldots,C_m$.
It is clear that $P_m$
is a $k$-ary tree with a staircase labeling.

By the above
procedure, we obtain a $k$-ary tree
$\phi(C_0,C_1,\ldots, C_{n-1})=P_{n-1}$,
labeled by $C_0, C_1, \ldots, C_{n-1}$.
It can be checked that the maps $\varphi$ and $\phi$
are inverses of each other.
 This completes the proof. \qed

\begin{figure}[h,t]
\setlength{\unitlength}{0.8mm}
\begin{center}
\begin{picture}(190,85)
\put(10,65){\circle*{1.5}}\put(8.5,59){$\emptyset$}
\put(15,65){\vector(1,0){15}}\put(19,67){$C_1$}
\put(40,73){\circle*{1.5}}\put(50,58){\circle*{1.5}}
\put(40,73){\line(2,-3){10}}\put(43,72){$\emptyset$}
\put(47,52){$(\textbf{2})$} \put(60,65){\vector(1,0){15}}
\put(64,67){$C_2$}\put(85,80){\circle*{1.5}}\put(95,65){\circle*{1.5}}
\put(85,80){\line(2,-3){10}}\put(85,50){\circle*{1.5}}
\put(85,50){\line(2,3){10}}\put(88,79){$\emptyset$}
\put(97,64){$(\textbf{2})$}\put(88,47){$(\textbf{2},\textbf{1})$}
\put(111,65){\vector(1,0){15}}\put(115,67){$C_3$}
\put(150,80){\circle*{1.5}}\put(160,65){\circle*{1.5}}
\put(150,50){\circle*{1.5}}\put(140,65){\circle*{1.5}}
\put(150,80){\line(2,-3){10}}\put(150,50){\line(2,3){10}}
\put(140,65){\line(2,3){10}}\put(153,78){$\emptyset$}
\put(162,64){$(\textbf{2})$}\put(152,48){$(\textbf{2},\textbf{1})$}
\put(132,58){$(\textbf{1},2,2)$} \put(15,25){\vector(1,0){15}}
\put(55,25){\circle*{1.5}}\put(35,24){$(\textbf{1},2,2)$}
\put(65,40){\circle*{1.5}}\put(75,25){\circle*{1.5}}
\put(65,10){\circle*{1.5}}\put(45,10){\circle*{1.5}}
\put(55,25){\line(2,3){10}}\put(65,40){\line(2,-3){10}}
\put(65,10){\line(2,3){10}}\put(45,10){\line(2,3){10}}
\put(67.5,38){$\emptyset$}\put(77,24){$(\textbf{2})$}
\put(32,2){$(\textbf{1},2,2,\textbf{1})$}
\put(60,2){$(\textbf{2},\textbf{1})$}\put(19,27){$C_4$}\put(90,25){\vector(1,0){15}}
\put(130,25){\circle*{1.5}}\put(110,24){$(\textbf{1},2,2)$}
\put(140,40){\circle*{1.5}}\put(150,25){\circle*{1.5}}
\put(140,10){\circle*{1.5}}\put(120,10){\circle*{1.5}}
\put(130,25){\line(2,3){10}}\put(140,40){\line(2,-3){10}}
\put(140,10){\line(2,3){10}}\put(120,10){\line(2,3){10}}
\put(142.5,38){$\emptyset$}\put(153,23){$(\textbf{2})$}
\put(106,2){$(\textbf{1},2,2,\textbf{1})$}
\put(135,2){$(\textbf{2},\textbf{1})$}\put(160,10){\circle*{1.5}}
\put(150,25){\line(2,-3){10}}
\put(94,27){$C_5$}\put(151,2){$(\textbf{2},\textbf{2},1,1,2)$}

\end{picture}\caption{An illustration of the bijection
$\phi$. }\label{p}
\end{center}
\end{figure}
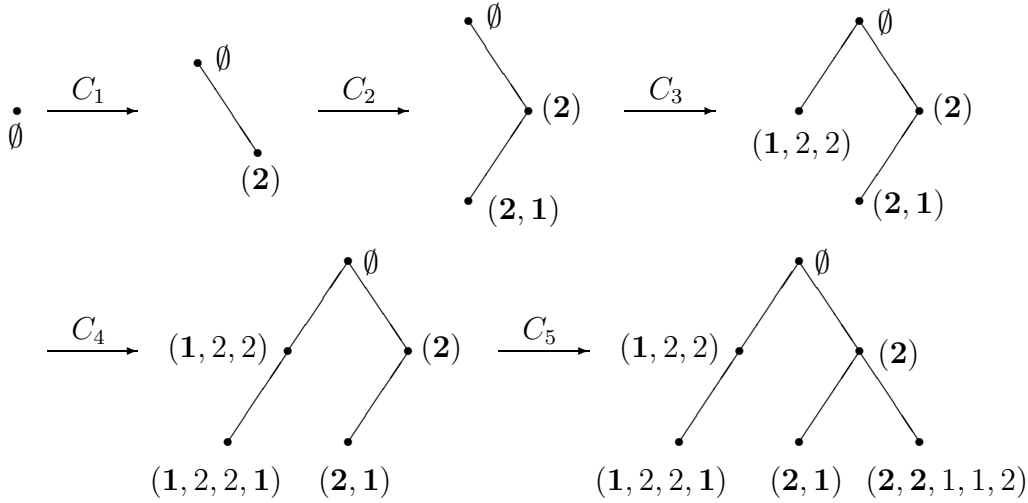

In particular,  for $k=2$, the proof of Theorem \ref{th}
reduces to a combinatorial proof of Han's formula (\ref{1})
 for binary trees. Figure \ref{p}
 gives an illustration of the bijection
 $\phi$ for $n=6$, $k=2$ and
\[(C_0,C_1,\ldots, C_5)=(\emptyset,(2),(2,1),(1,2,2),
(1,2,2,1),(2,2,1,1,2))\in S(6,2).\]

\section{A combinatorial interpretation of (\ref{2})}

In this section, we shall
apply the bijection $\phi$ constructed  in the
previous
section to give a combinatorial interpretation
of formula  (\ref{2}). To
this end, we need to reformulate
 (\ref{2}) in terms of complete binary
trees.

Clearly, one can add $n+1$ leaves
to a binary tree with $n$ vertices
to form a complete binary tree  with $2n+1$ vertices.
Moreover, a vertex
$u$ with hook length $h_u$ in a binary tree
becomes an internal
vertex with hook length
$2h_u+1$ in the corresponding complete binary tree.
Denote by $B^c_{2n+1}$ the set of complete binary trees
with $2n+1$ vertices.
Then (\ref{2}) is equivalent to
the  relation
\begin{align}\label{C}
\sum_{T\in B^c_{2n+1}} \frac{1}{\prod_{h\in
 \mathcal{H}(T)}h
2^{h-1}}=\frac{1}{2^n(2n+1)!}.
\end{align}

In fact, our combinatorial interpretation of
\eqref{C} is based on the following form
\begin{equation}\label{3}
\sum_{T\in B_{2n+1}^c}
\frac{(2n+1)!2^{1+2+\cdots+(2n+1)}}{\prod_{u\in \mathcal
{H}(T)}h2^h} =\frac{2^{1+2+\cdots +2n}}{2^n}.
\end{equation}

\noindent
{\it Combinatorial Proof of (\ref{3}).}
By the argument in the proof  of Lemma \ref{bb}, we see that the left-hand side of
(\ref{3}) is equal to the number of staircase labelings of complete
binary trees with $2n+1$ vertices. Let $S'(2n+1,2)$
the set of sequences in  $S(2n+1,2)$ corresponding to
complete binary trees with staircase labelings
under the bijection $\phi$. By
the construction of $\phi$,
we shall give an explanation of the fact that
\begin{equation}\label{7}
|S'(2n+1,2)|=\frac{1}{2^n}|S(2n+1,2)|.
\end{equation}
Since $|S(2n+1,2)|=2^{1+2+\cdots +2n}$,  we are led to a combinatorial
proof of  (\ref{3}).

It remains to prove (\ref{7}). To this end, we shall
construct a sequence of subsets $M_0, M_1, \ldots,
M_n$ such that
\[S(2n+1,2)=M_0\supset M_1 \supset\cdots \supset
M_n=S'(2n+1,2),\]
and for $1\leq i\leq n$,
\[|M_i|=\frac{1}{2}|M_{i-1}|.\]

Let us begin with the definition of
 the subset $M_1$ of $M_0$. Let $(C_0,C_1,\ldots,
C_{2n})$ be a sequence in  $M_0$, and
let $T$ be the corresponding binary tree under the
bijection $\phi$. If both subtrees of the
 root of $T$ (labeled with $C_0$) have an odd number of
 vertices, then we choose this sequence $(C_0, C_1,
 \ldots, C_{2n})$ to be in $M_1$.
By the construction of
$\phi$, it can be easily seen
that $(C_0,C_1,\ldots,
C_{2n})$ belongs to $M_1$ if and only
if there is an odd number of $1$'s among
$s_1,s_2,\ldots,s_{2n}$, where $s_i$ is
the first entry  of $C_i$. Since for any set
of an even number of elements, there are as many
subsets with an odd number of elements as subsets
with an even number of elements, we deduce that
\begin{equation}\label{m1}
|M_1|=\frac{1}{2}|M_0|.
\end{equation}

Similarly, we can define the subset $M_2$ of  $M_1$.
 Let $(C_0,C_1,\ldots,
C_{2n})$ be a sequence in $M_1$,
and let $T$ be the corresponding
binary tree under the
bijection $\phi$. Suppose that the vertices of $T$ are
$v_0, v_1, \ldots, v_{2n}$ and a vertex $v_i$ is labeled
by $C_i$. Clearly, $v_{0}$ is the root of $T$ labeled by $C_0$.
Suppose that $v_j$ is an non-root internal vertex with
$j$ being minimum. It is clear that $v_j$ is a child of $v_0$.
 If both subtrees of  $v_{j}$ have an odd number of
 vertices, then we choose this sequence $(C_0, C_1,
 \ldots, C_{2n})$ to be in $M_2$.
Using the above argument for (\ref{m1}) , we obtain that
\[|M_2|=\frac{1}{2}|M_1|.\]

In general, we can define the subset  $M_{j+1}$ of $M_{j}$
 for $j\geq 1$.
 Let $(C_0,C_1,\ldots,
C_{2n})$ be a sequence in $M_j$,
and let $T$ be the corresponding
binary tree under the
bijection $\phi$. Suppose that the vertices of $T$ are
$v_0, v_1, \ldots, v_{2n}$ such that a vertex $v_i$
is labeled
by $C_i$. Let $v_{t_0},v_{t_1},v_{t_2},\ldots$ be the
internal vertices  of $T$ such that the indices
are arranged in increasing order, that is,
 $t_0<t_1<t_2<\cdots$.
 If both subtrees of  $v_{t_j}$ have an odd number of
 vertices, then this sequence $(C_0, C_1,
 \ldots, C_{2n})$ is defined to be in $M_{j+1}$.
By the above reasoning, we see that
\[|M_{j+1}|=\frac{1}{2}|M_j|.\]

Once the subset $M_n$
has been determined, we get a binary tree
with a staircase labeling such that both subtrees
of any internal vertex have an odd number of vertices.
In other words, we obtain a complete binary tree
 with a staircase labeling. It follows that
 $M_n=S'(2n+1,2)$. This completes the proof. \qed

\vspace{.2cm} \noindent{\bf Acknowledgments.} This work was
supported by the 973 Project, the PCSIRT Project of the Ministry
of Education, and the National Science Foundation of China.



\begin{thebibliography}{99}

\bibitem {Che1}
W.Y.C. Chen, O.X.Q. Gao and P.L. Guo, Hook length formulas for trees
by  Han's expansion, Electron. J. Combin. 16 (1) (2009), R62.


\bibitem {Gel}
I.M. Gessel and  S. Seo, A Refinement of Cayley's formula for trees,
Electron. J. Combin. 11 (2) (2006), R27.


\bibitem {Han1}

G.-N. Han, Discovering new hook length formulas by expansion
technique, Electron. J. Combin. 15 (1) (2008), R133.

\bibitem {Han2}
G.-N Han, New hook length formulas for binary trees, Combinatorica
30 (2) (2010), 253-256.

\bibitem {Pos}
A. Postnikov, Permutohedra, associahedra, and beyond, Int. Math.
Res. Notices (2009), 1026-1106.

\bibitem {Sag}
 B.E. Sagan, Probabilistic proofs of hook length formulas involving
 trees, S\'em. Lothar. Combin., Special issue dedicated to the memory
  of Pierre Leroux, 61A (2009), Art. B61Aa.

\bibitem {Seo}
S. Seo, A combinatorial proof of Postnikov's identity and a
generalized enumeration of labeled trees, Electron. J. Combin. 11
(2) (2006), N3.

\bibitem {Yang}
 L.L.M. Yang, Generalizations of Han's hook length identities,
 arXiv:math.CO/0805. 0109.

\end{thebibliography}
\end{document}